\documentclass{au} 
\usepackage{eepic,epic}

\usepackage{a4wide}

\usepackage{amssymb,latexsym} 
\usepackage{mathrsfs}

\newcommand{\forget}[1]{}

\newcommand{\wg}[1]{\stackrel{(#1)}{=}}

\newcommand{\KA}{{\it KA}}

\newcommand{\Co}{{\mathbf{C}}}
\newcommand{\Loco}{\mathop{\mathrm{Loc_o}}}
\newcommand{\sSup}{\mathrm{sSup}}
\newcommand{\sub}{\subseteq}
\newcommand{\Sym}{\mathop{\mathrm{Sym}}}
\newcommand{\SymA}{\Sym(A)}
\newcommand{\Aut}{\mathop{\mathsf{Aut}}}
\newcommand{\sInv}{\mathop{\mathsf{sInv}}}
\let\sinv\sInv
\let\aut\Aut
\newcommand{\inv}{{\it inv}} 
\newcommand{\pot}{{\mathscr{P}}}
\newcommand{\M}{{\ul M}}
\newcommand{\RelmA}{\mathrm{Rel}^{(m)}(A)}
\newcommand{\RelA}{\mathrm{Rel}(A)}
\newcommand{\Rel}{\mathrm{Rel}}
\newcommand{\dom}{\mathop{\mathrm{dom}}}
\newcommand{\ima}{\mathop{\mathrm{im}}}
\newcommand{\arity}{\mathop{\mathrm{arity}}}
\newcommand{\Th}{\mathop{\mathit{Th}}}
\newcommand{\id}{\mathrm{id}}
\newcommand{\idA}{\id_A}
\newcommand{\cH}{\mathcal{H}} 
\newcommand{\cL}{\mathcal{L}} 
\newcommand{\cT}{\mathcal{T}} 
\newcommand{\cK}{\mathcal{K}}

\newcommand{\lab}[1]{\label{#1}%
}

\theoremstyle{plain}
\newtheorem{theo}{Theorem}[section]
\newtheorem{lemma}[theo]{Lemma}

\theoremstyle{definition}
\newtheorem{defi}[theo]{Definition}
\newtheorem{exa}[theo]{Example} 
\newtheorem{infd}[theo]{Informal Discussion}
\newtheorem{rema}[theo]{Remark}

\newcommand{\beao}{\begin{eqnarray*}}
\newcommand{\eeao}{\end{eqnarray*}}
\newcommand{\zwtext}[1]{}

\newcommand{\bea}{\begin{eqnarray}}
\newcommand{\eea}{\end{eqnarray}}
\newcommand{\be}{\begin{equation}}
\newcommand{\ee}{\end{equation}}

\newcommand{\eod}{\mbox{}\hfill\framebox[2mm]{\rule{0mm}{1pt}} \bigskip }

\newcommand{\eop}{\end{proof}}
\newcommand{\prf}{\begin{proof}}

\newcommand{\ul}[1]{\underline{#1}}

\newcommand{\pto}{\ \circ\!\!\!\to}
\renewcommand{\leq}{\leqslant}

\begin{document}
\title{Automorphisms  and strongly invariant relations}

\author{Ferdinand B\"orner} 
\author{Martin Goldstern} 
\thanks{The second author is grateful to the department of mathematics at 
Rutgers university for their hospitality during a visit in Fall 2002}
\author{Saharon Shelah}  %

\thanks{The third author is supported by NSF and the Israel Science
Foundation.  Publication 822.}
\date{2003-09-09} 
\begin{abstract}
We investigate characterizations of the Galois connection
        ${\rm sInv}$--${\rm Aut}$ between sets of finitary relations on a
        base set $A$ and their automorphisms. In particular, for
        $A=\omega_1$, we construct a countable set $R$ of relations that is
        closed
        under all invariant operations on relations and under arbitrary
        intersections, but is not closed under ${\rm sInv}{\rm Aut}$.
        
        Our structure $(A,R)$ has an  $\omega$-categorical first order
        theory.  A  higher order definable well-order makes it rigid,
        but any reduct to a finite language is homogeneous.
\end{abstract}

\maketitle

\section{Introduction} 

Our main question is easy to formulate. Let $R$ be a set of finitary 
relations on a nonempty base set $A$, and let $\Aut R$ denote 
the set of all automorphisms of the structure 
$(A\,;\,(\varrho)_{\varrho\in R})$. Conversely, if 
$G$ is a set of permutations on $A$, then $\sInv G$ denotes the 
the set of all relations $\sigma$ on $A$ such that all 
permutations in $G$ are automorphisms of $\sigma$. (Formal Definitions 
follow in the next section.) The question is: 
How can we characterize the relation sets  of the form 
$\sInv\Aut R$~? 

Of course, the operator $\sInv\Aut$ is a closure operator, and the 
operator pair $\sInv$--$\Aut$ forms a Galois connection between sets 
of relations on $A$ and sets of permutations on~$A$.   We can reformulate our 
problem as ``Which sets  $R$ of relations 
are {\em Galois-closed}, i.e., satisfy $R=\sInv\Aut R$?''
 or:  ``Describe the closure operator $\sInv\Aut$ {\em internally}'', i.e.,
without explicit reference to permutations.

Probably the first one who investigated this question in a systematic way 
was Marc Krasner. Influenced by the Galois connection 
between permutation groups and field extensions he tried to 
`generalize the notion of a field' \cite{Kra36}. Instead of the action 
of permutations on field elements, he considered the more complex action 
on relations. For finite base sets $A$ he described the closed sets 
of relations with the help of some operations on relations. A 
{\it logical operation on relations} is an operation, definable by 
a formula of the first order logic. (For details see the next section.) 
We call a set of relations a 
{\it Krasner algebra} if it 
is closed under all logical operations. For finite $A$, the Galois 
closed sets 
of relations are exactly the Krasner algebras. (At this point we 
remark that our notation differs from Krasner's original notation.) 

It is easy to extend this characterization to countable base sets~$A$: 
In this case the Galois closed relation sets are exactly those 
Krasner algebras 
that are additionally closed under arbitrary intersections 
($\bigcap$--closed Krasner algebras). 

But this is no longer true for the general case of  uncountable 
sets~$A$. For this case there exists a characterization by 
R. P\"oschel \cite{Poe80,Poe84} with the help of additional operations 
of uncountable arity. But the use of such operations is not 
very satisfying. Therefore we continue to look for better results. 

One reason for the existence of $\bigcap$--closed Krasner algebras 
that are not Galois closed is the fact that first order logic 
is   simply ``too weak''  to distinguish between sets of 
different infinite 
cardinalities. Consequently, it is a natural idea to replace the 
logical operations by a stronger class of operations. An $n$--ary 
operation $F$ on relations is called {\it invariant}, if the following 
identity holds for all permutations $g$ and all relations 
$\varrho_1,\varrho_n$ (with appropriate arities) on~$A$: 
$$ F(g[\varrho_1],\ldots,g[\varrho_n])= 
g\left[F(\varrho_1,\ldots,\varrho_n)\right] $$

Clearly, every Galois closed set of relations is $\bigcap$--closed 
and closed under all invariant operations. But it was unknown whether 
the converse is also true. The problem is: {\it Does there exist 
a set of relations that is $\bigcap$--closed and closed under all 
invariant operations, but not Galois closed for $\sInv$--$\Aut$?} 
(\cite[Problem 2.5.2]{Boe99}.) 

Surprisingly, the answer to this question is {\it yes}! In 
the main part of our article, section~3, we give a model theoretical 
construction of such a set of relations on a base set $A$ of cardinality 
$\omega_1$.

Finally, in section~4, we give a characterization of the Galois closed 
relation sets with the help of additional invariant infinitary 
operations. In contrast to P\"oschels characterization, we 
restrict these infinite arities to be countable. 
 Section~3 shows that we cannot restrict the arities to be finite, so
 this seems to be the best possible result.

\section{Preliminaries.} 

\subsection*{Notation} 

Throughout, let $A$ denote a nonempty base set. Write  $\omega$ for 
the set of all natural numbers (and at the same time the first infinite 
ordinal). An $m$--ary relation on $A$ 
is a subset of $A^m$, the set of all $m$--ary relations is denoted by  
$\RelmA$, and $\RelA:=\bigcup_{1\leq m\in\omega}\RelmA$ is the set of all 
finitary relations. If $R\sub\RelA$, then $R^{(m)}:=R\cap\RelmA$ 
. 
We do not distinguish between relations and predicates, 
therefore $\ul{a}\in\varrho$ and $\varrho(\ul{a})$ have the same meaning. 
The 
set of all permutations on $A$ is denoted by $\SymA$. For $g\in\SymA$ 
and $\ul{a}=(a_1,\ldots,a_m)\in A^m$ we put 
$$ g(\ul{a}) := (g(a_1),\ldots,g(a_m)), $$ 
and for $\varrho\sub A^m$ we write 
$$ g[\varrho] := \{g(\ul{a})\mid \ul{a}\in\varrho \}. $$

Let $g\in\SymA$ and $\varrho\in\RelA$. We say that $g$ is an 
{\it automorphism of $\varrho$}, or that $g$ {\it strongly preserves 
$\varrho$}, or that $\varrho$ is a {\it strongly invariant relation for  
$g$}, if 
$$ g[\varrho]=\varrho. $$
This is equivalent to $g[\varrho]\sub\varrho$ and 
$g^{-1}[\varrho]\sub\varrho$. 

For $R\sub\RelA$ and $G\sub\SymA$ we define operators 
$\Aut : \pot\RelA\to\pot\SymA$ and $\sInv : \pot\SymA\to\pot\RelA$:   
\beao 
\Aut R & := & \{ g\in \SymA \mid g[\varrho]=\varrho 
\mbox{ for all } \varrho\in R \} \\
\sInv G & := & \{ \varrho \in \RelA \mid g[\varrho] = \varrho \mbox{ for all } g\in\ G \} 
\eeao 
(For a set $X$, $\pot X$ denotes the set of all subsets of~$X$.) 

The operator pair $\sInv$--$\Aut$ forms a {\it Galois connection} between 
sets of permutations and sets of relations on $A$, i.e. the following 
conditions are satisfied: 
\begin{itemize}
\item $R_1\sub R_2\ \Rightarrow \ \Aut R_2 \sub \Aut R_1$ and 
$G_1 \sub G_2\ \Rightarrow \ \sInv G_2 \sub \sInv G_1$ 
\item $R\sub \sInv\Aut R$ and $G\sub\Aut\sInv G$.  
\end{itemize}
Consequently, the operators $$\sInv\Aut : \pot\RelA\to\pot\RelA\mbox{ and }
\Aut\sInv : \pot\SymA\to\pot\SymA$$ are closure operators. The sets of 
relations and the sets of permutations which are closed under these 
closure operators are called {\it Galois closed} (with respect to 
the Galois connection $\sInv$--$\Aut$). {\it Characterizing a Galois 
connection} means to describe the Galois closed sets without referring to 
the connection itself. 

In our article we want to find and discuss characterizations of our 
Galois connection $\sInv$--$\Aut$. There exist many similar Galois 
connections between 
sets of relations and sets of different kinds of functions, and they turned 
out to be useful especially for the investigation of finite mathematical 
structures. As a general source, we refer to \cite{PoeK79} and the 
list of references given there. Here we are interested in characterizations 
for infinite base sets~$A$. 

The main tool for the description of the closed sets of relations are 
{\it operations on relations}. These operations are of the form 
$$ F : \Rel^{(m_1)}(A)\times \ldots \times \Rel^{(m_n)}(A) \to \RelmA, $$ 
with $0\leq n\in\omega$. A set $R\sub\RelA$ is {\it closed under $F$} if 
$F(\varrho_1,\ldots,\varrho_n)\in R$ for all $\varrho_i\in R^{(m_i)}$, 
$1\leq i\leq n$. 

Special operations on relations are the {\it logical operations} which can 
be defined with the help of first order formulas. More exactly: Let 
$\varphi(P_1,\ldots,P_n;x_1,\ldots,x_m)$ be a formula with predicate 
symbols $P_i$ (of arity $m_i$), where all free  variables are in  
$\{x_j\mid 1\leq j\leq m\}$. We define 
$$ L_\varphi(\varrho_1,\ldots,\varrho_n) := 
\{ (a_1,\ldots,a_m)\in A^m\mid 
\varphi_A(\varrho_1,\ldots,\varrho_n, a_1,\ldots,a_m) \},  $$ 
where $\varphi_A(\varrho_1,\ldots,\varrho_n, a_1,\ldots,a_m)$ means that  
$\varphi$ holds in the structure 
$\langle A\,;\, \varrho_1,\ldots,\varrho_n\rangle$ for the 
evaluation $x_j:=a_j$, $(1\leq j\leq m)$. 

Examples of logical operations are the Boolean operations  
intersection $\cap$ and complementation $\Co$, 
defined by the formulas $P_1(x_1,\ldots,x_m)\wedge P_2(x_1,\ldots,x_m)$ 
and $\neg P(x_1,\ldots,x_m)$.

\subsection*{Properties of Galois closed relation sets} 

\begin{defi}\lab{defKA} 
A set $R\sub\RelA$ is called a {\it Krasner algebra (KA) on $A$}, 
if $R$ is closed 
under all logical operations.\footnote{An older notation is 
{\it Krasner algebra of second kind}, see e.g. \cite{PoeK79}.}   

If $Q\sub\RelA$, then $\langle Q\rangle_{\KA}$ denotes the {\it 
Krasner algebra generated by $Q$}, i.e. the least set of  relations 
on $A$ that contains $Q$ and is closed under all logical operations.  

A set $R$ of relations is called {\it $\bigcap$--closed}, 
if
 \begin{enumerate}
\item $A^m\in R$ for 
all $m\in\omega\setminus\{0\}$
\item  $R$ is closed under arbitrary 
intersections
 \end{enumerate}
 i.e. for all $m$ and all $Q\sub R^{(m)}$ we have 
$\bigcap Q \in R^{(m)}$. (Here we put $\bigcap\emptyset=A^m$.) 

If $Q\sub\RelA$, then $\langle Q\rangle_{\KA,\bigcap}$ denotes the 
least $\bigcap$--closed Krasner algebra containing~$Q$. 
\eod\end{defi}
Clearly, $\langle \ul{\mbox{\ \ }}\rangle_{\KA}$ 
and $\langle \ul{\mbox{\ \ }}\rangle_{\KA,\bigcap}$ are closure operators 
with $\langle Q\rangle_{\KA}\sub \langle Q\rangle_{\KA,\bigcap}$ 
for all $Q\sub\RelA$. 

If a set $R$ of relations is $\bigcap$--closed and closed under 
complementation $\Co$ (e.g. if $R=\langle R\rangle_{\KA,\bigcap}$), 
then it is also closed under arbitrary unions. 

The next Lemma gives the obvious connection between these 
notions and the Galois closed relation sets. 
For a proof we refer e.g. to \cite{Boe99}. 

\begin{lemma}\lab{lemKA}
If $R\sub\RelA$ is Galois closed ($R=\sInv\Aut R$), then $R$ is a 
$\bigcap$--closed Krasner algebra, $R=\langle R\rangle_{\KA,\bigcap}$. 
Consequently, for all $Q\sub\RelA$ we have
$\langle Q\rangle_{\KA,\bigcap}\sub\sInv\Aut Q$.    
\eod
\end{lemma} 
(Galois closed sets of relations are sometimes called 
{\it Krasner clones}. So every Krasner clone is a Krasner algebra, 
but not vice versa.) 

\begin{defi}\lab{defGamma}
Let $R\sub\RelA$ and $\ul{a}\in A^m$. We define: 
$$ \Gamma_R(\ul{a}) := 
\bigcap \{ \varrho \in R^{(m)} \mid \ul{a}\in\varrho \}\ \ \ \ \ \eod $$ 
\end{defi} 

We collect some properties of $\Gamma$. 
\begin{lemma}\lab{lemGamma} 
Let $R_1,R_2,R\sub\RelA$, $G\sub\SymA$ and $\ul{a}\in A^m$. 
Then the following hold.  
    \begin{enumerate}
    \item 
$R_1\sub R_2 \Rightarrow \Gamma_{R_2}(\ul{a})\sub\Gamma_{R_1}(\ul{a})$ 
    \item $\ul{a}\in\Gamma_R(\ul{a})$ and $\Gamma_R(\ul{a})\sub\varrho$ 
for all $\varrho\in R^{(m)}$ with $\ul{a}\in\varrho$. Moreover,  
$\varrho=\bigcup_{\ul{a}\in\varrho}\Gamma_R(\ul{a})$ for all $\varrho\in R$.
    \item If $R$ is $\bigcap$--closed, 
then $\Gamma_R(\ul{a})\in R$.  
\item If $R$ is closed under complementation, 
then $\{\Gamma_R(\ul{a})\mid \ul{a}\in A^m\}$ is a partition of $A^m$ 
and the relation $\ul{a}\sim_R\ul{b}:\iff \ul{a}\in\Gamma_R(\ul{b})$ is 
an equivalence relation. In this case, $R^{(m)}$ is an atomic Boolean 
algebra. 

  Note that   $\ul{a}\sim_R\ul{b}$ iff there is 
no relation $\varrho\in R$ separating $\ul{a} $ from   $\ul{b} $.

\item If $R_1$ and $R_2$ are $\bigcap$--closed and closed under 
complementation, then $R_1=R_2$ if and only if $\Gamma_{R_1}(\ul{a})=\Gamma_{R_2}(\ul{a})$ 
for all $m$ and all $\ul{a}\in A^m$. 
\item $\Gamma_{\sInv G}(\ul{a})=\{ g(\ul{a}) 
\mid g\in \langle G\rangle_{group}\}$, 
where $ \langle G\rangle_{group}$ is the subgroup of $\SymA$, 
generated by~$G$. 
    \end{enumerate}
\end{lemma}
\prf (1)--(5) are direct consequences of the Definitions. 
For (6), we first note that $\{ g(\ul{a}) 
\mid g\in \langle G\rangle_{group}\}$ contains $\ul{a}$ and 
is strongly invariant for all $g\in G$. Therefore 
$\Gamma_{\sInv G}(\ul{a})\sub\{ g(\ul{a}) 
\mid g\in \langle G\rangle_{group}\}$. On the other hand, 
$\sInv G$ is $\bigcap$--closed, therefore 
$\ul{a}\in\Gamma_{\sInv G}(\ul{a})\in\sInv G$, and every relation with these 
properties must contain all $g(\ul{a})$ with $g\in\langle G\rangle_{group}$. 
Consequently also $\Gamma_{\sInv G}(\ul{a})\supseteq\{ g(\ul{a}) 
\mid g\in \langle G\rangle_{group}\}$.
\eop

\begin{lemma}\lab{lemchauto}
Let $R\sub\RelA$ be $\bigcap$--closed and closed under complementation. 
Then $R=\sInv\Aut R$ if and  only if for all $m$ and all 
$\ul{a},\ul{b}\in A^m$ with $\ul{a}\sim_R \ul{b}$ there exists 
an automorphism $g\in\Aut R$ with $\ul{b}=g(\ul{a})$.  
\end{lemma}
\prf 
We have $R=\sInv\Aut R$ iff $\Gamma_R(\ul{a})=\Gamma_{\sInv\Aut R}(\ul{a})$ 
for all $\ul{a}$. Because of \ref{lemGamma}(6), this is equivalent to
$\Gamma_R(\ul{a})=\{ g(\ul{a})\mid g\in\Aut R\}$.  
\eop

\subsection*{Partial automorphisms} The last Lemma can be used to 
find characterizations in some special cases. 

\begin{defi}\lab{defpauto} 
A partial automorphism $f$ of a relation set $Q\sub\RelA$ (or 
of the structure $\ul{A}=( A ; (\sigma)_{\sigma\in Q})$)
with {\it domain} $\dom f=A_1\sub A$ and {\it image} $\ima f=A_2\sub A$ 
is a bijective function $f : A_1 \to A_2$, such that for all 
$\sigma\in Q$, $m = \arity(\sigma)$ and all $a_1,\ldots,a_m\in\dom f$, we have: $\sigma(a_1,\ldots,a_m)\Leftrightarrow\sigma(f(a_1),\ldots,f(a_m))$.

A set $Q\sub\RelA$ (or  the structure 
$\ul{A}=( A ; (\sigma)_{\sigma\in Q})$) is said to be 
{\it homogeneous}, if every finite partial automorphism can be extended to 
an automorphism of~$Q$.   
\eod\end{defi} 

\begin{lemma}\lab{lemhomog}
If $Q\sub\Rel( A)$ is a homogeneous  set of relations, then 
$\langle Q\rangle_{\KA,\bigcap}=\sInv\Aut Q$.  
\end{lemma}
\prf 
First note that  $\sInv \Aut Q$ is homogeneous and 
$ \langle Q\rangle_{\KA,\bigcap} \subseteq \sInv \Aut Q$,  
so also  $\langle Q\rangle_{\KA,\bigcap}$ is homogeneous.  So wlog $Q= 
 \langle Q\rangle_{\KA,\bigcap}$. 

Let $\ul{a}=(a_1,\ldots,a_m)$ and  
$\ul{b}=(b_1,\ldots,b_m)$.  We will use \ref{lemchauto}. So,
assuming $\ul a \sim_Q \ul b$, we have to find 
an automorphism $g$ with $g(\ul a) = \ul b$. 
Note that $\ul a \sim_Q \ul b$ in particular implies that $a_i=a_j$ iff $b_i=b_j$, for any $i,j$. So the map $f:=
\{(a_i, b_i)\mid i=1,\ldots, m\}$ is a finite
1-1 map.   As no relation in $Q$ separates $\ul a $ from $\ul b$, $f$ 
is even a partial automorphism for~$Q$. 

Because of the homogeneity of $Q$, $f$ can 
be extended to an automorphism $g\in\Aut Q$. Therefore $\ul{b}=g(\ul{a})$ 
for some $g\in\Aut Q$.  
\eop

Relation sets of the form $\sInv\Aut Q$ are homogeneous. Hence, 
a $\bigcap$--closed Krasner algebra is Galois closed if and only if it 
is homogeneous. 

\subsection*{The Galois closed permutation sets} 
We want to have a short look at the other side of our Galois connection. 
The characterization of the Galois closed permutation sets is 
well known (\cite{Jon68}) and provides no difficulties. We need an additional 
closure operator $\Loco : \pot\SymA\to\pot\SymA$: 
$$ \Loco G := \{ f\in\SymA \mid (\forall m\in\omega\setminus\{0\}) 
(\forall \ul{a}\in A^m ) (\exists g\in G) f(\ul{a})=g(\ul{a}) \} $$ 
\begin{theo}\lab{theocharperm}
A set $G\sub\SymA$ is Galois closed ($G=\Aut\sInv G$) if and only 
if $G=\langle G\rangle_{group}$ and $\Loco G=G$.  \eod 
\end{theo}
For the proof we refer to \cite{Jon68}. The 
operator $\Loco$ is a {\it topological closure operator}, 
multiplication and inversion of permutations are continuous with 
respect to the underlying topology. Therefore, the Galois closed 
automorphism sets are characterized as certain topological 
groups. For a more detailed discussion 
we refer to \cite[4.1]{Hod97}. 

\subsection*{A first characterization of the Galois closed relation sets}    
In \cite{Poe80} and \cite{Poe84}, R. P\"oschel characterized the 
closed sets of relations with the help of infinitary operations. 
Let $I$ be an arbitrary index set, let $m,m_i\in\omega\setminus\{0\}$ 
($i\in I)$. 
For an 
$I$-tuple $(\varrho_i)_{i\in I}$ of relations with $\varrho_i\in\Rel^{(m_i)}(A)$ 
the {\it strong superposition} with parameters  
$\ul{a}\in A^m$, $\ul{b}_i\in A^{m_i}$ is defined 
as follows: 
$$ \sSup_{\ul{a},(\ul{b}_i)_{i\in I}}(\varrho_i)_{i\in I} := 
\{ g(\ul{a}) \mid g\in\Sym(A) \mbox{ and } g(\ul{b}_i)\in\varrho_i 
\mbox{ for all } i\in I\} $$ 
A set $R\sub\RelA$ is {\it closed under strong superposition}, if 
$\sSup_{\ul{a},(\ul{b}_i)_{i\in I}}(\varrho_i)_{i\in I}\in R$ whenever $\varrho_i\in R$ 
for $i\in I$. 
\begin{theo}\lab{theopch}
Let $R\sub\RelA$ be  $\bigcap$--closed and closed under $\Co$. Then 
$R=\sInv\Aut R$ if and only if $R$ is closed under strong superposition. 
\end{theo}
\prf 
If $f\in\SymA$, then the definition of  
$\sSup_{\ul{a},(\ul{b}_i)_{i\in I}}$ implies 
$$ \sSup_{\ul{a},(\ul{b}_i)_{i\in I}}(f[\varrho_i])_{i\in I}=
f\left[\sSup_{\ul{a},(\ul{b}_i)_{i\in I}}(\varrho_i)_{i\in I}\right].  $$ 
Therefore every automorphism of  $\{\varrho_i \mid i\in I\}$ is an automorphism 
of $\sSup_{\ul{a},(\ul{b}_i)_{i\in I}}(\varrho_i)_{i\in I}$. Consequently, every Galois closed 
set of relations is closed under strong superposition. 

If $R$ is closed under strong superposition, then we can choose 
$I$ and $(\ul{b}_i)_{i\in I}$ such that all finite sequences with elements of $A$ 
occur among the $\ul{b}_i$. (This is possible with $|I|=|A|$.) 
Then $R$ contains the relation 
$\sSup_{\ul{a},(\ul{b}_i)_{i\in I}}(\,\Gamma_R(\ul{b}_i)\,)_{i\in I}$ and 
consequently $\Gamma_R(\ul{a})\sub  
\sSup_{\ul{a},(\ul{b}_i)_{i\in I}}(\,\Gamma_R(\ul{b}_i)\,)_{i\in I}$. 
This implies for every  $\ul{b}\in\Gamma_R(\ul{a})$ the existence 
of a permutation $g\in\SymA$ with $\ul{b}=g(\ul{a})$ and with 
$g(\ul{c})\in\Gamma_R(\ul{c})$ for all $n$ and all $\ul{c}\in A^n$. 
This $g$ is an automorphism of $R$, therefore \ref{lemchauto} 
implies that $R$ is Galois closed.  
\eop 

As seen in the proof, we can restrict the arities of the strong 
superpositions to $|I|=|A|$. Nevertheless, to be closed under strong 
superposition 
is a very strong condition. It immediately implies the existence of 
the necessary automorphisms in the sense of \ref{lemchauto}. Therefore, 
we continue to find better characterizations. 

\subsection*{A Characterization for countable base set $A$} 
For finite base set $A$, the Galois closed relation sets are exactly the 
Krasner algebras (\cite{Kra36,Kra50, Kra86,PoeK79}). This result can be 
extended to the countable case: 
\begin{theo}\lab{theocountA}
Let $A$ be a countable or finite set and $R\sub\RelA$. Then 
$R=\sInv\Aut R$ if and only if $R$ is a $\bigcap$--closed Krasner algebra, 
$R=\langle R\rangle_{\KA,\bigcap}$. 
Therefore $\langle Q\rangle_{\KA,\bigcap}=\sInv\Aut Q$ for all $Q\sub\RelA$. 
\end{theo} 
\prf
For the proof we 
refer to \cite[3.3.6.(v)]{Boe89} or to \cite[2.4.4.(i)]{Boe99}. 
One direction is provided by Lemma~\ref{lemKA}. For the other direction 
we can use a back \& forth construction to obtain the automorphisms  
that are necessary to apply Lemma~\ref{lemchauto}. 
\eop 

The next example shows that this characterization cannot be extended 
to uncountable sets.

\begin{exa}\label{ex1} Consider the following three countable structures: 
\begin{enumerate}
\item $({\mathbb Q}; {<})$ (the rational numbers with the linear order). 
\item The full countable bipartite graph: $(A\cup B; \varrho )$, where $A$
and $B$ are disjoint countable sets, and $\varrho = (A\times B) \cup
(B\times A)$.
\item The countable random graph. (See e.g. \cite[6.4.4]{Hod97}.) 
\end{enumerate}
Each of these structures $\M = (M;\varrho)$ has the following properties:
\begin{itemize}
\item[(a)] $Th(\M)$, the first order theory of $\M$, is $\omega$-categorical. 
\item[(b)] All unary first order formulas $\varphi(x)$ are equivalent
(mod $Th(\M)$) to $x=x$ or to $x\not=x$, i.e., the only subsets of
$\M$ that are first order definable without parameters are the empty
set and the whole model. 
\item[(c)] For any uncountable cardinal $\kappa$ there is a model
 $\M_\kappa$ of cardinality $\kappa$ such that the set
 $$\varrho^*:=\{ x: \mbox{The set $\{y: \varrho(x,y)\}$ is countable}\}$$
is neither empty nor the full model.
\end{itemize}

In each of these models $\M_\kappa$, the set $R$ of first order definable
relations (without parameters)
 is clearly a Krasner algebra and  is trivially closed under
$\bigcap$ (since,  by Ryll-Nardzewski's theorem, for any $k$ there 
are only finitely many $k$-ary relations in $R$). 

But in each model $\M_\kappa$ the set $\varrho^*$ is a (higher order)
definable subset of $M_\kappa$, hence $\varrho^*\in \sInv\Aut(R)\setminus R$.

This shows that $R$ is not Galois closed. 
\end{exa} 

\subsection*{Invariant operations} 
In the last example, the logical operations, together with arbitrary 
intersections, are too weak to provide the closure under $\sInv\Aut$. 
In particular, with logical operations it is not possible 
to distinguish between sets of different infinite cardinality. 
So the next idea to obtain a  
characterization is to replace the logical operations by a family 
of stronger operations on relations.  
\begin{defi}\lab{defIO} 
An operation 
$F : \Rel^{(m_1)}(A)\times \ldots \times \Rel^{(m_n)}(A) \to \RelmA$ is 
called {\it invariant}, if for all $g\in\SymA$ and all 
$\varrho_i\in\Rel^{m_i}(A)$  the following identity holds: 
$$ F(g[\varrho_1],\ldots,g[\varrho_n])=
g\left[F(\varrho_1,\ldots,\varrho_n)\right] $$

If $Q\sub\RelA$, then $\langle Q\rangle_{\inv}$ denotes the closure 
of $Q$ under all invariant operations, and $\langle Q\rangle_{inv,\bigcap}$ 
denotes the least set of relations that is closed under all invariant 
operations, $\bigcap$--closed and contains the set~$Q$. 
\eod\end{defi}
(The notations ``logical operations'' and ``invariant operations'' 
are adopted from \cite{Jon91}.) The operators 
$Q\mapsto\langle Q\rangle_{\inv}$ and $Q\mapsto\langle Q\rangle_{\inv,\bigcap}$ 
are closure operators, and we have 
$\langle Q\rangle_{\inv}\sub\langle Q\rangle_{\inv,\bigcap}$ for all 
$Q\sub\RelA$. 

We collect some easy properties of the invariant 
operations. 
\begin{lemma}\lab{leminvop}
\begin{enumerate}
\item Every logical operation is invariant. If $A$ 
is finite, then every invariant operation is logical. 
\item The invariant operations form a clone, i.e.\ the superposition of 
invariant operations is again an invariant operation. 
\item If $F$ is invariant and $\varrho_i\in\RelA$ ($1\leq i\leq n$), then 
$$\Aut\{\varrho_1,\ldots,\varrho_n\}\sub 
\Aut\{F(\varrho_1,\ldots,\varrho_n)\}.$$   
\item $ R \subseteq \langle R\rangle_{\inv}\subseteq \langle R\rangle_{\inv,\bigcap} \subseteq \sInv\Aut R$ for all~$R$. 
So, if $R=\sInv\Aut R$, then $R$ is $\bigcap$--closed and 
closed under all invariant operations, $R=\langle R\rangle_{\inv,\bigcap}$.
 
\end{enumerate}
\end{lemma}
\prf For (1) and (2) we refer to \cite{Jon91}. (3) is a direct consequence of 
\ref{defIO} and (4) is a direct consequence of (3). 
\eop 

The next Lemma shows that invariant operations are sufficient, 
if we have only 
finitely many relations. 
\begin{lemma}\lab{leminvopalg}
Let $Q\sub\RelA$ be a finite set. Then $\sInv\Aut Q=\langle Q\rangle_{\inv}$.  
\end{lemma}
\prf Let $Q=\{\varrho_1,\ldots,\varrho_n\}$ and let $\sigma\in\sInv\Aut R$. 
We define an invariant operation $F$ with 
$F(\varrho_1,\ldots,\varrho_n)=\sigma$:  
$$ F(\sigma_1,\ldots,\sigma_n):= \bigcup \{ g[\sigma] \mid 
g\in\SymA \mbox{ and } (\forall i=1\ldots n) g[\varrho_i]=\sigma_i \} $$ 
It is easy to verify that $F$ has the desired properties. \eop 

Let $\cH : \pot Z\to\pot Z$ be a closure operator on a set~$Z$. The 
{\it algebraic part of $\cH$} is the closure operator 
$$ \cH^{alg} :\pot Z\to\pot Z, \ X\mapsto\bigcup\{\cH X_0\mid 
X_0\sub X,\ X_0 \mbox{ finite}\}. $$
$\cH$ is {\it algebraic} if $\cH=\cH^{alg}$.

Lemma \ref{leminvopalg} shows 
that the closure operator  $Q\mapsto\langle Q\rangle_{\inv}$ 
is the {\it algebraic part} of $\sInv\Aut$. More exactly, for all 
$Q\sub\RelA$ we have
$$  \langle Q\rangle_{\inv} = \bigcup\{ \sInv\Aut Q_0\mid Q_0\sub Q 
\mbox{ and } Q_0 \mbox{ is finite}\} =(\sInv\Aut)^{alg} Q$$ 

So far, we have the closure operators $\langle\ul{\mbox{\ \ }}\rangle_{\KA}$, 
$\langle\ul{\mbox{\ \ }}\rangle_{\inv}$, $\langle\ul{\mbox{\ \ }}\rangle_{\KA,\bigcap}$, 
$\langle\ul{\mbox{\ \ }}\rangle_{\inv,\bigcap}$ and $\sInv\Aut$. 
The operators    $\langle\ul{\mbox{\ \ }}\rangle_{\KA}$ and  
$\langle\ul{\mbox{\ \ }}\rangle_{\inv}$ are algebraic, the other operators 
are not algebraic if $A$ is infinite. For all $Q\sub\RelA$ we have 
$$\langle Q\rangle_{\KA}\sub 
\begin{array}{c}
\langle Q\rangle_{\inv}\\[2mm]
\langle Q\rangle_{KA,\bigcap}
\end{array}
\sub\langle Q\rangle_{\inv,\bigcap}\sub\sInv\Aut Q.$$ 

For finite base set $A$ all these closure operators coincide. 
For countable $A$, the operators  
$\langle\ul{\mbox{\ \ }}\rangle_{\KA}$, 
$\langle\ul{\mbox{\ \ }}\rangle_{\inv}$ 
and  $\langle\ul{\mbox{\ \ }}\rangle_{\KA,\bigcap}$ are pairwise 
distinct -- see e.g. \cite[Theorem 4]{NemA88}, and 
$\langle\ul{\mbox{\ \ }}\rangle_{\KA,\bigcap}=
\langle\ul{\mbox{\ \ }}\rangle_{\inv,\bigcap}=\sInv\Aut$. If $A$ is 
uncountable, then also $\langle\ul{\mbox{\ \ }}\rangle_{\KA,\bigcap}\neq 
\langle\ul{\mbox{\ \ }}\rangle_{\inv,\bigcap}$, as a consequence of 
the examples in~\ref{ex1}.


All these properties now 
lead to the {\it conjecture}, that the Galois closed sets of relations 
are exactly the sets of relations that are $\bigcap$--closed 
and closed under all invariant operations. In order to verify this 
conjecture, we have to answer one question: 

{\it Does there exist a set $R$ of relations that is $\bigcap$--closed 
and closed under all invariant relations, but is not Galois closed? } 

This question was formulated as an open problem e.g. 
in \cite[Problem 2.5.2]{Boe99}. 

Surprisingly, the question has a positive answer and therefore the 
conjecture above is false. In the next section we will give a model 
theoretic construction of a relation set $R$ with the mentioned properties. 

\section{A model theoretic construction}
In this section we consider relational models of the form 
$\ul{M}=(M;(\varrho_m)_{1\leq m\in \omega})$, where 
$\varrho_m\in\RelmA$ for all~$m$. So our language $\cL$ has exactly 
one relation symbol for every arity~$m$. (We use the predicate 
symbols also to denote the corresponding relations.) 

If $\ul{M}=(M;(\varrho_m)_{1\leq m\in \omega})$ is such a 
model, then 
$\ul{M}^{[m]}:=
(M;\varrho_1,\ldots,\varrho_m)$ denotes the reduct 
of $\ul{M}$ to the relations $\varrho_1,\ldots,\varrho_m$. 

Our construction is guided by the following Lemma.

 \begin{lemma}\lab{lemmain0} 
 We fix  a vocabulary of infinitely many relational symbols 
 $\{\varrho_m \mid m \in \omega\}$. 

 Let $\ul{A} = ( A\,;\,(\varrho_m)_{1\leq m\in\omega})$ be an 
 infinite model.   
 Let $\ul{A}^{[m]} = (A; \varrho_1,\ldots, \varrho_m)$. 

 We assume that the following  hold: 
 \begin{enumerate}
 \item The theory $\Th(\ul{A})$ is $\omega$-categorical (i.e., has up to 
isomorphism exactly one countable model). 
 \item For all $m$, the reduct $\ul{A}^{[m]}$ is homogeneous in the sense of 
 \ref{defpauto}. 
 \item $\ul{A}$ is rigid, i.e. $\Aut \ul A =\{\id_A\}$.  
 \end{enumerate} 
 Then,  letting $R:=\langle \varrho_1,\varrho_2,\ldots\,\rangle_\KA$
 be the set of first order definable relations in $\ul A$, 
we have: 

$$\langle \varrho_1,\ldots, \varrho_m\rangle_{\KA} 
= \sInv \,\Aut \{ \varrho_1,\ldots, \varrho_m\}
$$
and $R= \bigcup_m \langle \varrho_1,\ldots, \varrho_m\rangle_{\KA} $, but
$$ 
R = \langle R\rangle_{\inv, \bigcap } \subsetneq \sInv\,\Aut R$$

 \end{lemma}

\begin{proof}
 First note that $R$ (as well as $\langle \varrho_1,\varrho_2,\ldots\varrho_m\rangle_\KA$) is closed under arbitrary
 intersections, since (by Ryll-Nardzewski's theorem)
there are only finitely many $k$-ary relations in $R$, for any~$k$.

We now show that $R$ is also closed under all invariant operations. 
We have
$$  
\langle \varrho_1,\varrho_2,\ldots,\varrho_m\rangle_\KA = 
 \langle \varrho_1,\varrho_2,\ldots,\varrho_m\rangle_{\KA,{\bigcap}}$$
Trivially,
$$ 
 \langle \varrho_1,\varrho_2,\ldots,\varrho_m\rangle_\KA
 \subseteq 
 \langle \varrho_1,\varrho_2,\ldots,\varrho_m\rangle_\inv
 \subseteq 
 \sinv\,\aut\{\varrho_1,\varrho_2,\ldots,\varrho_m\},$$
but by (2) and \ref{lemhomog} we have 
$$ \langle \varrho_1,\varrho_2,\ldots,\varrho_m\rangle_{\KA,{\bigcap}}
 =  \sinv\,\aut\{\varrho_1,\varrho_2,\ldots,\varrho_m\}, $$
so 
$$ 
\langle \varrho_1,\varrho_2,\ldots,\varrho_m\rangle_\KA = 
\langle \varrho_1,\varrho_2,\ldots,\varrho_m\rangle_\inv.$$

As both operators $ \langle\ul{\ \ }\rangle_{\KA}$ and 
 $ \langle\ul{\ \ }\rangle_{\inv}$  are algebraic,
 this yields 
$$ R = \langle \varrho_1,\varrho_2,\ldots\,\rangle_\KA = 
 \langle \varrho_1,\varrho_2,\ldots\,\rangle_\inv$$
hence $R = \langle R\rangle_{\inv, {\bigcap}}$.

Clearly $R$ is countable.  But $\Aut R=\{\idA\}$ and so $\sInv\Aut
R=\RelA$, which is an uncountable set. Consequently $$\langle
R\rangle_{\inv,\bigcap}\neq \sInv\Aut R.$$

\end{proof}

It remains to prove that a model with properties (1)--(3) exists. 
Because of \ref{theocountA}, such a model cannot be countable. 
We start by defining the  logical theory $\cT$  that we want our
model to satisfy. First we 
define an appropriate notion of a {\it clause}. 

\begin{defi}\lab{defclause}
A {\it literal in the variables $x_0,\ldots,x_n$} ($n\in\omega$) 
is a formula of the form 
$$  \varrho_m(x_{i_1},\ldots,x_{i_m}) \mbox{ (unnegated) } \mbox{ or } 
\neg\varrho_m(x_{i_1},\ldots,x_{i_m}) \mbox{ (negated) }$$
such that $1\leq m\leq n+1$, $\{i_1,\ldots,i_m\}\sub \{0,\ldots,n\}$, 
the  $i_1,\ldots,i_m$ are pairwise distinct and $0\in\{i_1,\ldots,i_m\}$. 

A {\it clause in $x_0,\ldots,x_n$} is a conjunction $K$ of literals in 
$x_0,\ldots,x_n$, such that no literal will appear twice, and no 
literal appears in negated and unnegated form. 

\eod\end{defi}  
Please note that there are only finitely many clauses in $x_0,\ldots,x_n$. 
The variable $x_0$ plays a special role --- it has to appear in every 
literal. 

Now we formulate our theory $\cT$: 
\begin{defi}\lab{deftheory}
$\cT$ consists of (the universal closures of) the
  following formulas: Firstly, for all $1\leq m\in\omega$ 
we have:
\begin{enumerate}
\item[(T1)] 
$\varrho_m(x_1,\ldots,x_m) \rightarrow 
\bigwedge_{1\leq i<j\leq m}x_i\neq x_j $
\end{enumerate}
Secondly, for all $n\in\omega$ and all clauses $K=K(x_0,\ldots,x_n)$ 
in $x_0,\ldots,x_{n}$ 
we take the formula: 
\begin{enumerate}
\item[(T2)] 
$ \bigwedge_{1\leq i<j\leq n} x_i\neq x_j \rightarrow 
(\exists x_0)\, K(x_0,\ldots,x_{n})  $   
\end{enumerate}
\eod
\end{defi} 

\begin{infd}\label{informal} \sl 
We will see below that the theory  $\cT$  is complete and
$\omega$-categorical. Our aim is to construct an uncountable model
$\ul M$  of
this theory on the base set $\omega_1$ in which the well-order 
$(\omega_1,{<})$  is definable (by a formula in higher order logic).    To
help us achieve this aim, we use the following ``recommendation'': 
\begin{quote} 
$\varrho(x_1,x_2,\ldots, x_n) $ {\bf should} hold iff $x_1 < x_2 <\cdots < x_n$
\end{quote}
However, this is just a recommendation, not a law.  In order 
to also get homogeneity of the restricted models $\M^{[m]}$,
we allow our model to disobey this recommendation, if there is a good
reason for it.   A good reason can be the desire to satisfy an axiom 
of our theory, or to extend a partial  automorphism. 

To keep track of the cases where the recommendation is not followed, 
we construct an auxiliary function $h: M \to \omega$, and we will
demand the following ``law'', which is a relaxed version of
 the ``recommendation'':
\begin{quote} 
For all sufficiently long tuples 
$(x_1,\ldots, x_n)$: \\
\ \qquad\qquad
$\varrho(x_1,x_2,\ldots, x_n) $ {\bf must hold} if $x_1 < x_2 <\cdots < x_n$,
and {\bf must not hold} otherwise
\end{quote}
Here, ``sufficiently long'' is defined as:  $n > \max (h(x_1), \ldots,
h(x_n))$. 

Thus, whenever we violate our recommendation at a tuple $(x_1,\ldots,
x_n)$, we will  define a sufficiently large value of $h$ at 
one of the points $x_1,\ldots, x_n$. 
\end{infd}

Before we investigate the theory $\cT$, we want to examine 
some technical definitions and lemmas. $\omega_1$ denotes the 
first uncountable ordinal, $\omega_1=\{\alpha\mid \alpha <\omega_1\}$. 
$\omega_1$ is well-ordered by~$<$. 
The universes of all our models will be subsets of $\omega_1$.

\begin{defi}\lab{defeckig}
Let $\ul{M}=( M;(\varrho_m)_{1\leq m\in\omega})$ be a model,  
let $h : M \pto \omega$ be a partial function, $n\in\omega\setminus\{0\}$ 
and let $\ul{a}=(a_1,\ldots,a_n)\in M^n$. We say that $\ul{a}$ is 
{\it a weak $n$-tuple for $h$} if 
 $\dom h \cap \{a_1,\ldots,a_n\}\neq\emptyset$ and 
 $$\max h (\ul{a}) := \max\{ h(a_i)\mid a_i\in\dom h\} < n.$$
All other $n$-tuples are called {\it strong for $h$}.

Now let  
$\ul{N}=( N;(\sigma_m)_{1\leq m\in\omega})$ and 
$\ul{M}=( M;(\varrho_m)_{1\leq m\in\omega})$ be models 
with $N\subset M\sub\omega_1$. 
Let $h : M \pto \omega$ be a partial function with $\dom h=M\setminus N$. 
We write $\ul{N} \sqsubset_h \ul{M}$ if 
\begin{enumerate}
\item $\ul{N}\leq \ul{M}$ ($\ul{N}$ is a submodel of $\ul{M}$), and 
\item for all $n$-tuples 
$(a_1,a_2,\ldots,a_n)\in M^n$ that are weak for $h$  the following 
condition holds: 
$$ \varrho_n(a_1,\ldots,a_n) \iff a_1<a_2<\ldots <a_n $$ 
(Here $<$ is the well-ordering on $\omega_1$.) 
\end{enumerate}
We write $\ul{N}\sqsubset\ul{M}$ if $\ul{N} \sqsubset_h \ul{M}$ for 
some partial function~$h$ with $\dom h=M\setminus N$. 
\eod 
\end{defi} 

\begin{lemma}\lab{lemsqtrans}
  \begin{enumerate}
  \item The relation $\sqsubset$ is transitive.
\item If $ \ul{M_0} \sqsubset \ul{M_1} \sqsubset \cdots $ is a chain of 
length $\omega$, and 
$\ul{M_\omega}$ is the directed 
union of this chain 
	(i.e., $M_\omega = \bigcup_{n\in \omega} M_n$, 
	and each $\M_n $ is also a submodel of $\ul{M_\omega}$), 
then $\ul{M_j}\sqsubset\ul{M_\omega}$ for all $j<\omega$. 
\item Similarly, if $(\ul{M_i})_{i\le\alpha}$ 
(where $\alpha$ is a limit ordinal) is a continuous chain of 
models (i.e., $M_i \le M_j$ for all $i\le j\le \alpha$, and 
for each limit $\delta \le \alpha$ we have $M_\delta = 
\bigcup_{i < \delta} M_i$), and  
$$\forall i < \alpha:  \ \ul{M_i}\sqsubset\ul{M_{i+1}},$$ 
then $\ul{M_j}\sqsubset\ul{M_\alpha}$ for 
all $j<\alpha$. 
  \end{enumerate}
\end{lemma} 
\prf \

\begin{enumerate}
\item Assuming $\ul{M_1}\sqsubset_{h_1}\ul{M_2}$ and
$\ul{M_2}\sqsubset_{h_2}\ul{M_3}$ for some $h_1 : M_2\setminus
M_1\to\omega$ and $h_2 : M_3\setminus M_2\to\omega$, we have to show
$\ul{M_1} \sqsubset \ul{M_3}$.
\\
Put $h := h_1\cup h_2 : M_3\setminus M_1 \to \omega$. If 
$\ul{a}\in (M_3^n\setminus M_1^n)$ is weak for $h$, then 
either $\ul{a}\in (M_2^n\setminus M_1^n)$ or one of the $a_i$ belongs 
to $M_3\setminus M_2$. In the first case, 
$\max h_2(\ul{a})=\max h(\ul{a})<n$, 
and $\ul{a}$ is weak for $h_2$. Therefore ($\ul{M_2}$ is a submodel of 
$\ul{M_3}$), 
$$\varrho_n^{3}(\ul{a})\iff\varrho_n^{2}(\ul{a}) \iff a_1<a_2<\ldots<a_n.$$
In the second case, $\max h_3(\ul{a}) \leq \max h(\ul{a}) <n$, therefore 
$\ul{a}$ is weak for $h_3$ and $\varrho_n^{3}(\ul{a})\iff a_1<\ldots<a_n$.  
\item[(2)] and (3) similar.
\forget{ \tiny
\item Clearly, all $\ul{M_i}$ are submodels of $\ul{N}$. We have 
$\ul{M_i}\sqsubset_{h_i}\ul{M_{i+1}}$ for some 
$h_i : M_{i+1}\setminus M_i\to\omega$. Put $h := \bigcup_{j\leq i<\alpha}h_i$.
Then $h : N\setminus M_j \to\omega$. Let $\ul{a}\in(N^n\setminus M_j^n)$ be 
weak for~$h$. Let $k'<\alpha$  denote the least ordinal with 
$\{a_1,\ldots,a_n\}\sub M_k$. Our $n$ is finite, therefore $k'$ is not a 
limit ordinal, i.e. $k'=k+1$, and at least one of the $a_i$ is in 
$M_{k+1}\setminus M_k$. As in the case before, we have 
$\max h_k(\ul{a})\leq\max h(\ul{a})<n$, so $\ul{a}$ is already weak 
for $h_k$. Therefore $\varrho_n^{k+1}(\ul{a})\iff a_1<\ldots <a_n$. 
$\ul{M_k}$ is a submodel of $\ul{N}$, so this result also holds for 
$\varrho_n^{\ul{N}}$.   
} 
\end{enumerate}
\eop

The next technical Lemma provides the basic step in our construction. 
\begin{lemma}\lab{lemconstr}
Let $\ul{M_{0}}=( M_0 ; (\varrho_m^{0})_{1\leq m\in\omega})$ 
be a countable model of our language $\cL$, such that $M_0\sub\omega_1$ and 
$\varrho_m(x_1,\ldots,x_m) \rightarrow 
\bigwedge_{1\leq i<j\leq m}x_i\neq x_j$ holds for all~$m$. Moreover, let 
$\pi_0 : M_0 \pto M_0$ be a partial (finite or infinite) 
automorphism of the reduct 
$\ul{M_0}^{[s]}$ for some $s\in\omega$ and let $\alpha\in\omega$. 

Then there exists a countable model 
$\ul{M_{\omega}}=
( M_\omega ; (\varrho_m^{\omega})_{1\leq m\in\omega})$ with 
$M_0\subset M_\omega\subset\omega_1$, $\alpha\in M_\omega$ and 
a total
automorphism $\pi_\omega : M_\omega\pto M_\omega$ of the 
$s$--reduct $\ul{M_\omega}^{[s]}$ 
such that the following hold: 
\begin{enumerate}
\item $\ul{M_0}\sqsubset \ul{M_\omega}$ 
\item $\ul{M_\omega}$ is a model of the theory $\cT$. 
\item $\pi_\omega$ extends $\pi_0$ 
\end{enumerate}
\end{lemma}
\prf
We construct $\ul{M_\omega}$ as the union of a chain 
of countable models $\ul{M_j}$ with $j\in\omega$. We will have 
$\ul{M_j}\sqsubset\ul{M_{j+1}}$ for all $j\in\omega$, 
and for every $j$ we will have a partial automorphism 
$\pi_j$ of $\ul{M_j}^{[s]}$ such that $\pi_{j+1}$ extends 
$\pi_j$, and $\dom(\pi_{j+1})\cap \ima(\pi_{j+1})\supseteq M_j$. 

We explain the step from $M_j$ to $M_{j+1}$. Let 
$M_j=\{a_i\mid i\in\omega\}$  be an enumeration of the elements 
of $M_j$. (The elements $a_i$ are not necessarily in the 
order, given by $<$ in $\omega_1$.) Let 
$\cK$ be the following set: 
$$ \cK := \{ (n,\ul{a},K)\mid n\in\omega, 
\ul{a}\in M_j^n \mbox{ and $K$ a 
clause in } x_0,\ldots,x_n \} $$ 
This set is countable. Let $(n_l,\ul{a}_l,K_l)_{l\in\omega}$ be an 
enumeration of this set. 

Let $B\sub\omega_1\setminus M_j$ be a countable set of ordinals, 
and let $B=\{b_k\mid k\in\omega\}$ be an enumeration of~$B$. 
(Again, this enumeration need not necessarily follow the well-order 
$<$ on $\omega_1$.) We  put $M_{j+1} := M_j\cup B$, and we have to 
define the relations $\varrho_m^{j+1}$ and the partial 
function $\pi_{j+1}$. Moreover, we must define a function $h : B\to \omega$, 
in order to establish the relation $\ul{M_j}\sqsubset_h\ul{M_{j+1}}$.  

For all $m$ and all $m$-tuples $\ul{a}\in M_j^m$ we define: 
$$ \varrho_m^{j+1}(\ul{a}) :\iff \varrho_m^{j}(\ul{a}) $$
This makes sure that $\ul{M_j}\leq\ul{M_{j+1}}$. Moreover, 
we define $\neg\varrho_m^{j+1}(c_1,\ldots,c_m)$ for all 
$c_1,\ldots,c_m\in M_{j+1}$ with $|\{c_1,\ldots,c_m\}|<m$. 

For every $k\in\omega$ we will conduct a special task, 
where we define the value $h(b_k)$, define a partial function 
$p_k : M_{j+1}\pto M_{j+1}$ such that $p_{k+1}$ always extends $p_k$. 
We start with $p_0 := \pi_j$, and finally we will have 
$\pi_{j+1}:=  \bigcup_{k}p_k$. Moreover, in every step we 
define the 
truth values of $\varrho_m^{j+1}(\ul{a})$ for some tuples 
$\ul{a}$.

We distinguish three cases of steps $k$, 
depending on whether $k\equiv 0,1,\mbox{ or }2 
\mod (3)$.    [A main point will be that the definitions in the 
various cases do not contradict each other.] 
\bigskip

{\it Step $k$ for $k=3l$:} If $l=0$, then let $p_k:=\pi_j$, otherwise 
$p_k:=p_{k-1}$ remains unchanged. 
 
\begin{figure}

\setlength{\unitlength}{0.0004in}
\begingroup\makeatletter\ifx\SetFigFont\undefined%
\gdef\SetFigFont#1#2#3#4#5{%
  \reset@font\fontsize{#1}{#2pt}%
  \fontfamily{#3}\fontseries{#4}\fontshape{#5}%
  \selectfont}%
\fi\endgroup%
{ \renewcommand{\dashlinestretch}{30}
\begin{picture}(5100,4737)(0,-10)
\path(1212,3237)(5562,3237)(5562,12)
	(1212,12)(1212,3237)
\path(1512,2637)(3162,2637)(3162,2337)
	(1512,2337)(1512,2637)
\path(2787,2712)(4737,4662)
\path(4673.360,4555.934)(4737.000,4662.000)(4630.934,4598.360)
\put(3387,2412){\makebox(0,0)[lb]{$ \ul{a}_l$}}
\put(4062,3612){\makebox(0,0)[lb]{$ K_l$}}
\put(4962,4587){\makebox(0,0)[lb]{$b_l $}}
\put(1812,2412){\makebox(0,0)[lb]{$ \cdots \ \cdots $}}
\put(5962,2412){\makebox(0,0)[lb]{$M_j$}}
\put(-2400,2412){\makebox(0,0)[lb]{ $k=3l$}}
\end{picture}
}
\end{figure}

Let $(n_l,\ul{a}_l,K_l)$ be the element of $\cK$ with index~$l$. 
We define $h(b_k):=n_l+1$. Now, for every literal which occurs in $K$, 
we define the truth values in such a way, that $K_l(b_k,\ul{a}_l)$ 
becomes true.   Thus, letting $\ul{a}_l = (a_l(1),\ldots, a_l(n_l))$, 
we define truth values for certain tuples
from the set $\{a_l(1),\ldots, a_l(n_l),b_k\}^{<\omega} \setminus  
\{a_l(1),\ldots, a_l(n_l)\})^{<\omega}$.

\forget{\tiny
More exactly, if $\varrho_m(x_{i_1},\ldots,x_{i_m})$ 
occurs in $K_l$ and $c_0:=b_k$ 
and $c_i:=a_i$ for $i\geq 1$, then we define  
$\varrho_m^{j+1}(c_{i_1},\ldots,c_{i_m})$, and if 
$\neg\varrho_m(x_{i_1},\ldots,x_{i_m})$ 
occurs in $K_l$, then we define  
$\neg\varrho_m^{j+1}(c_{i_1},\ldots,c_{i_m})$. 
}

The largest index $m$ of a literal which occurs in $K_l$ is $n_l+1$.
Therefore all these tuples $\ul{c}$ satisfy $\max h(\ul{c})\geq
h(b_k)\geq m$ and are strong for~$h$.
[Note that in no previous step have we committed ourselves to the truth value of $\varrho_j(\ul{c})$ for any tupel $\ul c$ in which $b_k$ appears.]

\bigskip

{\it Step $k$ for $k=3l+1$:} In this case we define $h(b_k):=s$, and we extend 
the partial function $p_{k-1}$. If $\dom p_{k-1}\supseteq M_j$, then we 
simply put $p_k:=p_{k-1}$ and we are done. 

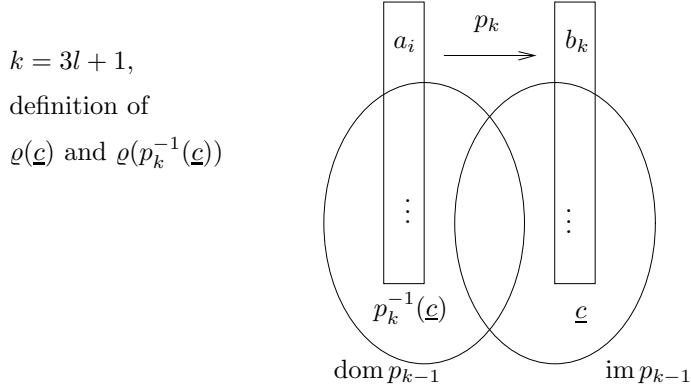
\begin{figure}

\setlength{\unitlength}{0.0004in}
\begingroup\makeatletter\ifx\SetFigFont\undefined%
\gdef\SetFigFont#1#2#3#4#5{%
  \reset@font\fontsize{#1}{#2pt}%
  \fontfamily{#3}\fontseries{#4}\fontshape{#5}%
  \selectfont}%
\fi\endgroup%
{\renewcommand{\dashlinestretch}{30}
\begin{picture}(4802,5214)(0,-10)
\put(3462,3012){\ellipse{1500}{2100}}
\put(2487,3012){\ellipse{1500}{2100}}
\path(12,3237)(4362,3237)(4362,12)
	(12,12)(12,3237)
\path(987,2937)(2862,5187)
\path(931,2984)(1043,2890)
\path(2808.225,5075.608)(2862.000,5187.000)(2762.131,5114.019)
\path(2337,4187)(2338,4189)(2341,4192)
        (2345,4198)(2351,4207)(2360,4219)
        (2371,4232)(2383,4248)(2398,4264)
        (2414,4281)(2432,4299)(2452,4316)
        (2474,4333)(2499,4350)(2527,4367)
        (2560,4383)(2597,4398)(2637,4412)
        (2675,4422)(2711,4430)(2744,4436)
        (2773,4441)(2797,4444)(2817,4446)
        (2833,4448)(2846,4449)(2858,4449)
        (2868,4449)(2879,4449)(2890,4449)
        (2903,4448)(2919,4446)(2938,4444)
        (2961,4441)(2988,4436)(3019,4430)
        (3053,4422)(3087,4412)(3123,4398)
        (3154,4383)(3182,4367)(3204,4350)
        (3224,4333)(3240,4316)(3254,4299)
        (3266,4281)(3277,4264)(3286,4248)
        (3293,4232)(3300,4219)(3304,4207)(3312,4187)

\path(3239.579,4287.275)(3312.000,4187.000)(3295.287,4309.559)
\put(3162,4487){\makebox(0,0)[lb]{$ p_{k{-}1}$}}
\put(987,2587){\makebox(0,0)[lb]{$ a_i$}}
\put(2937,5037){\makebox(0,0)[lb]{$b_k $}}
\put(1580,4062){\makebox(0,0)[lb]{$p_k $}}
\put(4587,2337){\makebox(0,0)[lb]{$ M_j$}}
\put(1012,1587){\makebox(0,0)[lb]{$ \dom p_{k{-}1}$}}
\put(3237,1587){\makebox(0,0)[lb]{$ \ima p_{k{-}1}$}}
\put(-2400,2412){\makebox(0,0)[lb]{ $k=3l+1$}}
\end{picture}
}
\end{figure}

If $M_j\setminus \dom p_{k-1}\neq\emptyset$, then let 
$i:=\min\{i\mid a_i\notin \dom p_{k-1}\}$. We  extend $p_{k-1}$ by  defining 
$p_k(a_i):=b_k$. Then, for all $m\leq s$ and all $c_1,\ldots,c_m\in\ima p_k$ 
such that $c_1,\ldots,c_m$ are pairwise distinct and 
$b_k\in\{c_1,\ldots,c_m\}$, we define  the truth value of 
$\varrho_m^{j+1}(c_1,\ldots,c_m)$. Write $\ul{c} $ for $(c_1,\ldots,c_m)$.
If the value of 
$ \varrho^{j+1}(p_k^{-1}(\ul{c}))$ is already known, 
in particular if $p_k^{-1}(\ul{c})\in M_j^m$, then we  put 
$$ \varrho_m^{j+1}(\ul{c}) :\iff 
\varrho^{j+1}(p_k^{-1}(\ul{c}).$$ 
If $ \varrho^{j+1}(p_k^{-1}(\ul{c}))$ is still not fixed 
(in particular, then the $p_k^{-1}(c_i)$ cannot all be in $M_j$), 
then we define both values, namely we put 

\begin{figure}

\setlength{\unitlength}{0.0007in}
\begingroup\makeatletter\ifx\SetFigFont\undefined%
\gdef\SetFigFont#1#2#3#4#5{%
  \reset@font\fontsize{#1}{#2pt}%
  \fontfamily{#3}\fontseries{#4}\fontshape{#5}%
  \selectfont}%
\fi\endgroup%
{\renewcommand{\dashlinestretch}{30}
\begin{picture}(2491,3084)(0,-10)
\put(1733,1200){\ellipse{1500}{2100}}
\put(758,1200){\ellipse{1500}{2100}}
\path(458,2850)(758,2850)(758,750)
	(458,750)(458,2850)
\path(1733,2850)(2033,2850)(2033,750)
	(1733,750)(1733,2850)
\path(908,2450)(1583,2450)
\path(1463.000,2420.000)(1583.000,2450.000)(1463.000,2480.000)
\put(83,0){\makebox(0,0)[lb]{$ \dom p_{k-1}$}}
\put(2108,0){\makebox(0,0)[lb]{$\ima p_{k-1}$}}
\put(1883,450){\makebox(0,0)[lb]{$ \ul c$}}
\put(383,450){\makebox(0,0)[lb]{$p_{k}^{-1}(\ul c) $}}
\put(608,1200){\makebox(0,0)[lb]{$ \vdots$}}
\put(1808,1125){\makebox(0,0)[lb]{$ \vdots$}}
\put(528,2475){\makebox(0,0)[lb]{$ a_i$}}
\put(1808,2475){\makebox(0,0)[lb]{$b_k $}}
\put(1133,2625){\makebox(0,0)[lb]{$ p_k$}}
\put(-2400,2312){\makebox(0,0)[lb]{ $k=3l+1$, }}
\put(-2400,2012){\makebox(0,0)[lb]{ definition of }}
\put(-2400,1612){\makebox(0,0)[lb]{ $\varrho(\ul{c})$ and $\varrho(p_k^{-1}(\ul{c}))$}}

\end{picture}
}
\end{figure}

$$ \varrho_m^{j+1}(\ul{c}) , 
\varrho_m^{j+1}(p_k^{-1}(\ul{c})):\iff 
p_k^{-1}(c_1)<\ldots <p_k^{-1}(c_m).$$ 
(Here $<$ denotes the well-order in $\omega_1$.) 

Because of the definition 
of $h(b_k)$, the $(c_1,\ldots,c_m)$ are always strong for~$h$ and hence are 
exempted from our ``recommendation'' \ref{informal}. 
 The other tuples, $p_k^{-1}(\ul{c})$, follow our recommendation anyway, whether or not they are $h$-strong. 

[As before, note that tuples $\ul c$ in which $b_k$ appears have never been 
considered in any previous step $k'< k$.]

\bigskip 

{\it Step $k$ for $k=3l+2$:} This step is similar to the previous
step, but this time we take care of $\ima p_k$ rather than $\dom p_k$,
or in other words: we reverse the roles of $p_k$ and $p_k^{-1}$.  We
leave the details to the reader. 

\forget{\tiny 
Again, we define $h(b_k):=s$, and we extend the partial function
$p_{k-1}$. If $\ima p_{k-1}\supseteq M_j$, then we simply put
$p_k:=p_{k-1}$ and we are done.

If $M_j\setminus \ima p_{k-1}\neq\emptyset$, then let 
$i:=\min\{i\mid a_i\notin \ima p_{k-1}\}$. We  extend $p_{k-1}$ by  defining 
$p_k(b_k):=a_i$. Then, for all $m\leq s$ and all $c_1,\ldots,c_m\in\dom p_k$ 
such that $c_1,\ldots,c_m$ are pairwise distinct and 
$b_k\in\{c_1,\ldots,c_m\}$, we define  the truth value of 
$\varrho_m^{j+1}(c_1,\ldots,c_m)$. If the value of 
$ \varrho^{j+1}(p_k(c_1),\ldots,p_k(c_m))$ is already known, then we put 
$$ \varrho_m^{j+1}(c_1,\ldots,c_m) :\iff 
\varrho_m^{j+1}(p_k(c_1),\ldots,p_k(c_m)).$$ 
If $ \varrho(p_k(c_1),\ldots,p_k(c_m))$ is still not fixed 
(in particular, then the $p_k(c_i)$ cannot all be in $M_j$), 
then we define both values, namely we put 
$$ \varrho_m^{j+1}(c_1,\ldots,c_m),
\varrho_m^{j+1}(p_k(c_1),\ldots,p_k(c_m)):\iff 
p_k(c_1)<\ldots <p_k(c_m).$$ 
(Again, $<$ denotes the well-order in $\omega_1$.) Because of the definition 
of $h(b_k)$, the $(c_1,\ldots,c_m)$ are always strong for~$h$. Again, the 
other tuples, $(p_k(c_1),\ldots,p_k(c_m))$ are correctly defined and cannot appear 
a second time in one of our steps. 
}
\bigskip
\bigskip

By induction we obtain from these steps a model, where the relations 
$\varrho_m^{j+1}$ are only partially defined. In order to finish the definition, 
we put  
$$ \varrho^{j+1}_m(\ul{c}) :\iff c_1<c_2<\ldots< c_m, $$
whenever the truth value of $ \varrho^{j+1}_m(\ul{c})$ has not been defined
during the inductive construction.

\relax From the construction, it is now clear that $\ul{M_j}\sqsubset\ul{M_{j+1}}$. 
Moreover, $\pi_{j+1}:=\bigcup_{k\in\omega} p_k$ is a partial automorphism 
of $\ul{M_{j+1}}^{[s]}$ that extends $\pi_i$ and satisfies 
$M_j\sub\dom\pi_{j+1}$ and $M_j\sub\ima\pi_{j+1}$. 
Moreover, all formulas in $\cT$ of the form \ref{deftheory}(T2)  are satisfied,
whenever $x_1,\ldots, x_m\in M_j$.

The countable set $B\subset\omega_1\setminus M_j$ was arbitrary. 
So, if the ordinal $\alpha$ is not in $M_0$, then we can assume that $\alpha\in B$, 
e.g. for $j=0$. Consequently we will have $\alpha\in M_\omega$ in the end.  

We form the directed union  $\ul{M_\omega}:=\bigcup_{j\in\omega}\ul{M_j}$ 
Because of \ref{lemsqtrans}(2), we have $\ul{M_0}\sqsubset\ul{M_\omega}$. 
Moreover, the union $\pi_\omega:=\bigcup_{j\in\omega}\pi_j$ is a 
bijective partial automorphism of $\ul{M_\omega}^{[s]}$, 
which is everywhere defined and surjective, i.e. 
it is an automorphism of $\ul{M_\omega}^{[s]}$  which extends $\pi_0$. 

Finally, if $\psi(x_1,\ldots,x_n)$ is a formula of $\cT$ and $a_1,\ldots,a_n\in M_\omega$, 
then there exists $j$ with $a_1,\ldots,a_n\in M_j$. Consequently $\psi(a_1,\ldots,a_n)$ 
holds in $\ul{M_{j+1}}$ and all other extensions of $\ul{M_j}$, in particular it 
is true 
in $\ul{M_\omega}$. Consequently 
$\ul{M_\omega}$ is a model of $\cT$. This finishes the proof.  
\eop

Now we collect some properties of our theory $\cT$. 
\begin{lemma}\lab{lemtheory} 
  \begin{enumerate}
  \item $\cT$ is consistent and has no finite models. 
  \item $\cT$ has the property of elimination of quantifiers. 
  \item $\cT$ is complete. 
  \item $\cT$ is $\omega$--categorical. 
  \item If $\ul{M},\ul{N}$ are models of $\cT$ and $\ul{N}$ is a submodel 
of $\ul{M}$, $\ul{N}\leq\ul{M}$, then $\ul{N}$ is an elementary submodel of 
$\ul{M}$, i.e. for every formula $\varphi(x_1,\ldots,x_n)$ and every 
$\ul{a}\in N^n$ holds $\varphi(\ul{a})$ in $\ul{N}$ iff it holds in $\ul{M}$.  
  \end{enumerate}
\end{lemma}
\prf
\begin{enumerate}
\item Easy. The consistency of $\cT$ is a corollary of \ref{lemconstr}, and the nonexistence 
of finite models is ensured  by the formulas \ref{deftheory}(T2).

\item In order to prove that for every formula $\psi(x_1,\ldots,x_n)$ there is a 
quantifier free formula $\varphi(x_1,\ldots,x_n)$ with 
$\cT\models (\varphi\leftrightarrow \psi)$, it is sufficient to prove this for 
formulas $\psi$ of the form $(\exists x_0) K(x_0,x_1,\ldots,x_n)$, where 
$K$ is a clause in $x_0,\ldots,x_n$. 

For any equivalence relation $\theta$ on $\{1,\ldots,n\}$ we put: 
$$ \mu_\theta := (\bigwedge_{(i,j)\in\theta}x_i=x_j)\wedge(\bigwedge_{(i,j)\notin\theta}x_i\neq x_j).$$
Let $E_n$ denote the set of all equivalence relations on $\{1,\ldots,n\}$. Then 
$$K\leftrightarrow \bigvee_{\theta\in E_n}(K\wedge\mu_\theta). $$ 
Therefore also 
$$(\exists x_0) K(x_0,x_1,\ldots,x_n) \leftrightarrow \bigvee_{\theta\in E_n}(\exists x_0)(K\wedge\mu_\theta).$$ 
It is sufficient to show that every formula $(\exists x_0)(K\wedge\mu_\theta)$ is equivalent 
to a quantifier free formula. If $\theta$ is not the equality relation, then we 
can replace any variable $x_j$  ($j\geq 1$) 
by $x_i$, where $i$ is a representative of the equivalence class of~$j$. If then 
a variable appears twice in an literal of $K$, then either the clause becomes false 
modulo $\cT$ (if the literal is unnegated), or the literal can be omitted modulo $\cT$ 
(if it is negated). In the end we obtain either formulas which are true (modulo $\cT$) 
or false (modulo $\cT$) or equivalent to a formula of the form 
$$  (\exists x_0)(K\wedge\bigwedge_{1\leq i<j\leq n}x_i\neq x_j).$$
$\cT$ contains the formula $\bigwedge_{1\leq i<j\leq n}x_i\neq x_j\rightarrow (\exists x_0)K$, 
therefore 
$$\cT\models (\bigwedge_{1\leq i<j\leq n}x_i\neq x_j\leftrightarrow (\exists x_0)(K\wedge\bigwedge_{1\leq i<j\leq n}x_i\neq x_j).$$
\item By (2), every closed formula is (mod $\cT$) equivalent to 
$\mathbf{true}$ or $\mathbf{false}$.  

\item Modulo $\cT$, there are only finitely many quantifier-free 
formulas in the variables $x_1,\ldots, x_n$, namely, 
Boolean combinations of atomic formulas
 $\varrho_m(x_{i_1}, \ldots, x_{i_m})$, for
 $i_1,\ldots, i_m\in \{1,\ldots, n\}$ and $m\le n$.
  (Note that formulas 
 $\varrho_m(x_{i_1}, \ldots, x_{i_m})$ for $m>n$ and $i_1,\ldots, i_m\le n$ are 
automatically false mod $\cT$, because of \ref{deftheory}(T1). 
\\
This implies $\omega$-categoricity, by Ryll-Nardzewski's theorem. 
 [Actually, we do not need $\omega$-categoricity itself for our construction, 
we only need the fact 
that there are only finitely many first order definable 
$k$-ary relations, for any $k$.] 

\forget{\tiny
The $\omega$--categoricity is a consequence of the fact that there are 
only finitely many {\it types} in the variables $x_1,\ldots,x_n$ (See e.g. \cite{Hod97}). 
Every type can be written as a single conjunction of literals 
$(\neg)\varrho_n(x_{i_1},\ldots,x_{i_n})$ with $n\leq m$,  
and there are only finitely many possible conjunctions of these finitely many literals.  
 
The completeness is a consequence of (2). (Every sentence is modulo $\cT$ equivalent to 
true or to false.)  }
\item This is a consequence from the fact that $\cT$ has elimination of quantifiers. 
\end{enumerate}
\eop

The results in \ref{lemtheory} make sure that the condition 
\ref{lemmain0}(1) is satisfied for every model of $\cT$. 
Now we construct a model $\ul{A}$, such that also the conditions 
\ref{lemmain0}(2) and (3) are satisfied.

We will obtain $\ul{A}$ as a directed union over an uncountable chain of models,  
$\ul{A}:=\bigcup_{i\in\omega_1}\ul{M_i}$, such that every  $\ul{M_i}$ is 
a model of $\cT$, $\ul{M_i}\sqsubset\ul{M_{i+1}}$ and $i\in M_{i+1}\subset\omega_1$  
for all $i\in\omega_1$. Because of \ref{lemtheory}(4), this is an elementary 
chain, therefore $\ul{A}$ is again a model of $\cT$. Because of $i\in M_{i+1}$ for 
all $i\in\omega_1$ and $M_i\sub\omega_1$, the carrier set $A=\bigcup_{i\in\omega_1}M_i$  
is  $A=\omega_1$.

In order to obtain a model with homogeneous $s$--reducts, we have to make sure 
that certain partial automorphisms can be extended to automorphisms. For this reason, 
we use a triply-indexed family $(\pi_{n,i,j})_{n\in\omega,i,j\in\omega_1, i\leq j}$ 
of partial automorphisms. 

First we explain, what the $\pi_{n,i,i}$ are.  
If $\ul{M}$ is a countable model, then there are only countable many pairs 
$(p,s)$ with $s\in\omega$ and $p$ a finite partial automorphism of the 
$s$--reduct $\ul{M}^{[s]}$. Therefore there exists an enumeration $(p_n,s_n)_{n\in\omega}$ 
of all these finite partial automorphisms with corresponding~$s$.  Now, for  $i\in\omega$ 
and 
$\ul{M}=\ul{M_i}$ we put $\pi_{n,i,i}:=p_n$, and $s_{n,i}:=s_n$. Therefore 
\begin{itemize}
\item[$(**)$]
$(\pi_{n,i,i})_{n\in\omega}$ is a list of all finite partial automorphisms 
of all possible reducts $\ul{M_i}^{[s]}  $. 
\end{itemize}
 The  $\pi_{n,i,j}$ with $i<j$ will be 
extensions of $\pi_{n,i,i}$. 

Now we explain how to construct the models  $\ul{M_j}$
and sequences of partial 
automorphisms  $(\pi_{n,i,j}: i \le j)$ by transfinite induction 
on $j \in \omega_1$). 
 This construction will use the usual `bookkeeping--argument' to
take
care of $\omega_1\times \omega$ many tasks in $\omega_1$ steps. 
Let $\omega_1=\bigcup_{n\in\omega,i\in\omega_1}C_{n,i}$ be a partition of $\omega_1$ 
into pairwise disjoint sets $C_{n,i}$, such that $|C_{n,i}|=\omega_1$ for all $(n,i)$ 
and $\min C_{n,i}\geq i$. 

If $j=0$, then let $\ul{M_0}$ be a countable model of $\cT$ with $M_0\sub\omega_1$. 
(The existence of such a model is clear from \ref{lemconstr}.) The $\pi_{n,0,0}$ are 
defined as in $(**)$.

If $j$ is a limit ordinal, then put $\ul{M_j}:=\bigcup_{i<j}\ul{M_i}$. 
(As a directed union of an elementary chain of models of $\cT$, this is again 
a model of $\cT$.) The $\pi_{n,j,j}$ are defined as in $(**)$, and  
$\pi_{n,i,j}:=\bigcup_{i\leq l< j}\pi_{n,i,l}$. 

For a successor ordinal $j+1$ we use Lemma~\ref{lemconstr}:
\begin{enumerate}
\item  We define $\ul M_{j+1}$ as follows.  Let $(n,i)$ be 
the pair with $j \in C_{n,i}$. According to the Lemma, there exists a model 
$\ul{M_{j+1}}$ with $j\in M_{j+1}\subset \omega_1$ and $\ul{M_j}\sqsubset\ul{M_{j+1}}$, 
and there exists an extension of $\pi_{n,i,j}$ to a partial automorphism 
$\bar \pi$ of $\ul{M_{j+1}}^{[s_{n,i}]}$ 
with $M_j\sub\dom\bar \pi$ and $M_j\sub \ima \bar \pi$. 
\item We let $\pi_{n,i,j+1}$ be the  partial automorphism $ \bar \pi$ 
from~(1). 
\item 
The $\pi_{n,j+1,j+1}$ are defined as in $(**)$, 
enumerating all finite partial automorphisms of reducts of $\ul{M_{j+1}}$. 
\item 
For all $(n,i)$ such that $j\notin C_{n,i}$, we put $\pi_{n,i,j+1}:=\pi_{n,i,j}$. 
\end{enumerate}

\bigskip

It is easy to verify by transfinite induction, that the $\pi_{(n,i,j)}$ are always 
partial automorphisms of $\ul{M_j}^{[s_{n,i}]}$, and that $\pi_{n,i,j}$ extends 
$\pi_{n,i,k}$ for all $k$ with $i\leq k<j$.

As mentioned above, we put 
$\ul{A}:=\bigcup_{i\in\omega_1}\ul{M_i}$. 

\begin{lemma}\lab{lemhomg2}
If $s\in\omega$, then every finite partial automorphism $\pi$ of $\ul{A}^{[s]}$ 
can be extended to an automorphism of $\ul{A}^{[s]}$. 
\end{lemma}
\prf
$\pi$ is finite, therefore there exists $i\in\omega$ with $\dom\pi\cup\ima\pi\sub M_i$. 
$\ul{M_i}$ is a submodel of $\ul{A}$, therefore $\pi$ is a finite partial automorphism of $\ul{M_i}^{[s]}$. 
Consequently, $\pi=\pi_{n,i,i}$ for some $n$ with $s=s_{n,i}$. We define 
$\pi':=\bigcup \{ \pi_{n,i,j} : {j\in \omega_1, j\geq i}\}$.
 $\pi'$ is a partial isomorphism of~$\ul{A } ^{[s]}$. 
By our construction, we have $M_j\sub\dom\pi'$ for all $j\in C_{n,i}$, i.e. 
$\dom\pi\supseteq \bigcup_{j\in C_{n,i}}M_j=\omega_1$. The same holds for 
$\ima \pi'$, therefore $\pi'$ is a total  automorphism  of~$\ul{A } ^{[s]}$.  
\eop

It remains to show that $\ul{A}$ has the property \ref{lemmain0}(3). 
\begin{lemma}\lab{lemrigid}
$\ul{A}=( A;(\varrho_m)_{m\in\omega\setminus\{0\}})$ 
is rigid, i.e. $\Aut\ul{A}=\{\idA\}$. 
\end{lemma}
\prf 
Let  $S$ be a countable subset of $A$, $x,y\in S$ and $h : A\setminus S\to\omega$ 
be a function. Then we define that 
$E(x,y,S,h)$ is true iff  for all $m$ and for all 
$\ul{a}=(a_1,\ldots,a_m)\in A^m\setminus S^m$ the following holds: 
\begin{quote}
If  $ (\varrho_m(\ul{a})\wedge \max h(\ul{a})<m \wedge 
(\exists i,j\in\{1,\ldots,m\}) (x=a_i \wedge y=a_j)  )$, \\
then $i<j$
\end{quote}
We claim  $x<y\iff (\exists S)(\exists h) E(x,y,S,h)$. 

\medskip
\noindent Proof of ``$\Rightarrow$'': Let $i\in\omega_1$ be the least  ordinal with $x,y\in
M_i$. 
Let $S:= M_i$. We have 
$\ul{M_i}\sqsubset\ul{A}$, therefore $\ul{M_i}\sqsubset_h\ul{A}$ for some 
$h : A\setminus S\to\omega$. If $\ul{a}\in A^m\setminus S^m$, is weak for $h$, 
then $\varrho_m(\ul{a})\iff a_1<a_2<\ldots<a_m$. Therefore, if $x=a_i$, $y=a_j$ and 
$x<y$, then $i<j$. 

\medskip
\noindent Proof of ``$\Leftarrow$'': Let $S$ be a countable subset of $A$ and let $h : A\setminus S\to \omega$ 
be a function such that $E(x,y,S,h)$. 

Since $\omega_1$ has uncountable cofinality, and $S$ is countable, there
must be some $i < \omega_1$ with $S \subseteq M_i$. 

So let $i$ be the least ordinal with $S\sub M_i$. Let $p\in A\setminus M_i$. (Then also 
$p\in A\setminus S$.) We have $\ul{M_i}\sqsubset_{h_i}\ul{A}$ for some  
$h_i : A\setminus M_i\to \omega$. Let $m:=\max\{h(p),h_i(p)\}+3$ and choose 
$z_1,\ldots,z_{m-3}\in S$, pairwise distinct and distinct from $x$ and~$y$. 
Let $\ul{a}=(a_1,\ldots,a_m)$ be the $m$-tuple consisting of the 
elements of $\{p,x,y,z_1,\ldots,z_{m-3}\}$ in the ordering according to $<$, 
i.e. $a_1<a_2<\ldots<a_m$.

Find $i$, $j$ such that 
$$ x=a_i, \ y=a_j.$$
 Thus, 
$$i<j \iff x< y. $$ 

First we note that 
$\max h_i(\ul{a})<m$, so $\ul{a}$ is weak for $h_i$.
  As $\ul{M_i}\sqsubset \ul{M}$, we must have $\varrho_m(\ul{a})$. 

But we also have $\max h(\ul{a})<m$, so $\varrho_m(\ul{a})$ implies
$i<j$.  So $x < y$.

\bigskip 

Let $\bar E(x,y) :\iff (\exists S)(\exists h) E(x,y,S,h)$.  If $\pi$
is an automorphism of $\ul{A}$ we write $\pi[h]$ for the map $h'$
satisfying $h'(\pi(x)) = h(x)$. 
Clearly  $E(x,y,S,h) \iff E(\pi(x),\pi(y), \pi[S], \pi[h])$, hence 
$$ \bar E(x,y) \iff \bar E(\pi(x), \pi(y)). $$

Consequently every automorphism of $\ul{A}$ has to preserve the well-ordering 
$<$ of $\omega_1$. But the only order automorphism of $<$ is $\idA$.   
\eop

In \ref{lemtheory}, \ref{lemhomg2} and \ref{lemrigid} we have verified all properties of $\ul{A}$, 
required in \ref{lemmain0}. So we can formulate our main result. 

\begin{theo}\lab{theomain}
The structure $\ul{A}=( A\,;\,(\varrho_m)_{m\in\omega})$ has the properties 
\ref{lemmain0}(1)--(3). Consequently, the set $R=\{\varrho_m\mid m\in\omega\}$ satisfies 
$$ \langle R\rangle_{\KA}=\langle R\rangle_{\inv,\bigcap}\neq\sInv\Aut R.$$ 
 \eod
\end{theo}

\begin{rema}
The structure $\ul{A}$ that we constructed 
in \ref{defclause}--\ref{lemrigid} 
has cardinality 
$\aleph_1$. A similar construction can be carried out 
to yield a model of any cardinality $\kappa$ with cofinality 
$\mathit{cf}(\kappa)\geq\omega_1$.  We leave the details 
to the reader. 
\end{rema} 


\section{A characterization with invariant operations of countable arity}
The results of the last section show that the closure under 
$\langle\ul{\mbox{\ \ }}\rangle_{\inv,\bigcap}$ is 
too weak to 
provide the closure under $\sInv\Aut$. So  the question remains, 
what we should add to obtain an appropriate characterization. 
Similar as in \ref{theopch} we will add operations with infinite arity. 
But in contrast to \ref{theopch}, we need only operations with 
countable arities. 

\begin{lemma}\lab{lemonerel}
Let $m\in\omega\setminus\{0\}$ and let $Q\sub\RelmA$ be closed under 
complementation. Then there exists a relation 
$\varrho\in\langle Q\rangle_{\KA,\bigcap}^{(2m)}$
with $\Aut Q=\Aut\{\varrho\}$. 
\end{lemma}
\prf 
$Q$ is closed under $\Co$, therefore the set 
$$M:=\{\Gamma_Q(\ul{a})\mid \ul{a}\in A^m\}$$  
is a partition of $A^m$. $M$ can be well-ordered, so let 
$(\gamma_i)_{i<\kappa}$ be a corresponding enumeration of the elements 
of~$M$. We put 
$$\varrho := \bigcup_{i\leq j<\kappa} \gamma_i\times\gamma_j $$ 
Then $\varrho\in\langle Q\rangle_{\KA,\bigcap}^{(2m)}$. 
This implies $\varrho\in\sInv\Aut Q$ and 
therefore $\Aut\{\varrho\}\supseteq\Aut Q$. 

Now let  $g\in\Aut\{\varrho\}$. If $\ul{a},\ul{b}\in\gamma_i$, then 
$(\ul{a},\ul{b})\in\varrho$ and $(\ul{b},\ul{a})\in\varrho$, 
hence $(g(\ul{a}),g(\ul{b}))\in\varrho$ and $(g(\ul{b}),g(\ul{a}))\in\varrho$.
The $\gamma_i$ are pairwise disjoint, therefore there are unique 
ordinals $l,k<\kappa$ with $g(\ul{a})\in\gamma_l$, $g(\ul{b})\in\gamma_k$. 
Now $(\ul{a},\ul{b})\in\varrho$ implies $l\leq k$ and 
$(g(\ul{a}),g(\ul{b}))\in\varrho$ implies $k\leq l$, i.e. $l=k$. 
Consequently all tuples in $\gamma_i$ are transformed by $g$ to 
tuples in $\gamma_k$. Therefore there exists a function 
$g_0 : \kappa\to\kappa$ with $g[\gamma_i]\sub\gamma_{g_0(i)}$. 

$g^{-1}$ is also an automorphism of $\varrho$, and it is easy to see that 
the corresponding function $g'_0 : \kappa \to \kappa$ has to be the inverse 
of $g_0$. Consequently, $g_0$ is a permutation on $\kappa$.  

The permutation $g_0$ preserves the well-order $<$ on $\kappa$, because of 
$$ i<j\Rightarrow\gamma_i\times\gamma_j\sub\varrho\Rightarrow 
\gamma_{g_0(i)}\times\gamma_{g_0(j)}\sub\varrho\Rightarrow g_0(i)<g_0(j). $$
But $\langle\kappa,<\rangle$ has only the trivial order automorphism, 
therefore $g_0=\id_\kappa$. 

We obtain $g[\gamma_i]\sub\gamma_i$ and (because $g^{-1}$ is also 
an automorphism) $g^{-1}[\gamma_i]\sub\gamma_i$, i.e. 
$g[\gamma_i]=\gamma_i$. But then, $g$ is an automorphism for all 
relations in $M$, and therefore also for all relations in~$Q$. 
This yields $\Aut\{\varrho\}\sub\Aut Q$, and  
this finishes the proof. 
\eop 

As a consequence of this Lemma, every possible 
automorphism group appears already as 
the automorphism group of an at most countable set of relations. 

\begin{lemma}\lab{lemvorletztes} 
For every set $R\sub\RelA$ there exists an at most countable set 
$R_0\sub\RelA$ with $\Aut R=\Aut R_0$. Moreover, if $R$ is a 
$\bigcap$--closed Krasner algebra, then we can choose $R_0\sub R$. 
\end{lemma} 
\prf 
If $R$ is not closed under $\Co$, then we put $R':=R\cup\{\Co\sigma\mid 
\sigma\in R\}$. Then $\Aut R= \Aut R'$, therefore we can assume 
w.l.o.g. that $R$ is closed under complementation. By \ref{lemonerel} 
there are relations $\varrho_m\in\Rel^{(2m)}(A)$ with 
$\Aut R^{(m)}=\Aut\{\varrho_m\}$. Consequently:  
$$ \Aut R = \bigcap_{1\leq m\in\omega} \Aut R^{(m)} = 
\bigcap_{1\leq m\in\omega} \Aut \{\varrho_m\} = \Aut \{\varrho_m\mid 1\leq m\in\omega\} $$
The second part in \ref{lemonerel} implies that the $\varrho_m$ 
can be chosen from the $\bigcap$--closed Krasner algebra, generated by~$R$. 
\eop

Now we define our additional operations. 
\begin{defi}\lab{defctblinv} 
An {\it invariant operation with countable arity} is an operation 
of the form 
$$ F : \prod_{1\leq i\in\omega}\Rel^{(m_i)}(A) \to \RelmA $$ 
($m_i\in\omega\setminus\{0\}$), such that for all 
$(\varrho_i)_{1\leq i\in\omega}\in \prod_{1\leq i\in\omega}\Rel^{(m_i)}(A)$ 
and all $g\in\SymA$ we have
$$ F(g[\varrho_i])_{1\leq i\in\omega}=
g\left[F(\varrho_i)_{1\leq i\in\omega}\right] $$ 
If $Q\sub\RelA$, then $\langle Q\rangle_{\omega-\inv}$ is the closure of 
$Q$ under all invariant operations with countable arity, and 
$\langle Q\rangle_{\omega-\inv,\bigcap}$ is the least set of relations 
which is closed under all invariant operations with countable arity and 
$\bigcap$--closed. 
\eod\end{defi} 
(Of course, $\langle \ul{\ \ }\rangle_{\omega -\inv}$ and 
$\langle \ul{\ \ }\rangle_{\omega -\inv,\bigcap}$ are closure operators.) 
Similar as in \ref{leminvop} and \ref{leminvopalg}, we can verify 
the following properties: 
\begin{lemma}\lab{lemom-inv}
\begin{enumerate}
\item If $R\sub\RelA$ is Galois closed, $R=\sInv\Aut R$, then 
$R$ is $\bigcap$--closed and closed under all invariant operations 
with countable arity, $R=\langle R\rangle_{\omega -\inv,\bigcap}$. 
\item If $Q\sub\RelA$ is countable or finite, then 
$\langle Q\rangle_{\omega -\inv}=\sInv\Aut Q$. 
\end{enumerate}
\eod
\end{lemma}

Now we can formulate our characterization of the Galois closed 
sets of relations. 

\begin{theo}\lab{theoletztes} 
Let $R$ be a $\bigcap$--closed Krasner algebra. Then $R$ is 
Galois closed, $R=\sInv\Aut R$, if and only if $\sInv\Aut R_0\sub R$ 
for every countable subset $R_0$ of~$R$. 

In particular, a set $R\sub\RelA$ is Galois closed 
if and only if it is $\bigcap$--closed and closed under all invariant 
operations with countable arity, $R=\langle R\rangle_{\omega -\inv,\bigcap}$. 
For all $Q\sub\RelA$ holds $\langle Q\rangle_{\omega -\inv,\bigcap}=\sInv\Aut Q$. 
\end{theo} 
\prf 
Clearly, $R_0\sub R$ and $R=\sInv\Aut R$ implies $\sInv\Aut R_0\sub\sInv\Aut R =R$ for every subset $R_0$ of~$R$.  Vice versa, if $\sInv\Aut R_0\sub R$ for 
all countable subsets, then we can choose the special subset $R_0\sub R$ 
 with $\Aut R=\Aut R_0$ of Lemma~\ref{lemvorletztes}. Then we obtain: 
$$ R\sub\sInv\Aut R = \sInv\Aut R_0 \sub R, $$  
i.e., $R$ is Galois closed. 

Then the other statements are consequences of Lemma \ref{lemom-inv}. 
\eop

\end{document}